\documentclass[11pt,twoside,a4paper]{amsart}
\usepackage{hyperref}

\usepackage[english]{babel}
\usepackage{amssymb}
\usepackage{amsthm}
\usepackage[all]{xy}
\usepackage{euscript}
\newcommand{\E}{\mathrm E}
\newcommand{\V}{\mathrm V}
\newcommand{\N}{\mathbb N}
\newcommand{\R}{\mathbb R}
\newcommand{\Z}{\mathbb Z}
\newcommand{\Hh}{\mathrm {H_0}}
\newcommand{\te}{\mathrm {\theta}}
\newcommand{\Uu}{\EuScript U}
\newcommand{\Vv}{\EuScript V}
\newcommand{\Ww}{\EuScript W}
\newcommand{\Xx}{\EuScript X}

\def\diam{\mathop{\rm diam}\nolimits}
\def\mod{\mathop{\rm mod}\nolimits}
\def\int{\mathop{\rm Int}\nolimits}
\def\cat{\mathop{\rm CAT(0)}\nolimits}

\def\andim{\mathop{\rm dim_{AN}}\nolimits}
\def\ldim{\mathop{\rm {\ell\hbox{-}dim}}\nolimits}
\def\lasdim{\mathop{\rm {\ell\hbox{-}asdim}}\nolimits}
\def\asdim{\mathop{\rm {asdim}}\nolimits}
\newcommand{\dist}{\operatorname{dist}}

\theoremstyle{plain}
\newtheorem{thm}{Theorem}
\newtheorem{sle}{Corollary}
\newtheorem{lmm}{Lemma}
\newtheorem{utv}{Proposition}
\newtheorem*{hyp}{Hypothesis}

\theoremstyle{definition}
\newtheorem{dfn}{Definition}

\theoremstyle{remark}
\newtheorem{zam}{Remark}
\newtheorem*{constr}{Construction}

 \newtheoremstyle{break}
   {9pt}
   {9pt}
   {\itshape}
   {}
   {\bfseries}
   {.}
   {\newline}
   {}

 \theoremstyle{break}

\begin{document}
\def\objectstyle{\scriptstyle}
\def\labelstyle{\scriptstyle}

\title[l-asdim of the fundamental group of a
graph-manifold]{Linearly-controlled asymptotic
dimension of the fundamental group of a
graph-manifold}
\author{Alexander Smirnov}
\subjclass{Primary 57M50, 55M10; Secondary 05C05, 20F69}

\begin{abstract}
We prove the estimate
$\lasdim \pi_1(M)\leq 7$
for the linearly controlled asymptotic dimension of the
fundamental group of any 3-dimensional graph-manifold
$M$.
As applications we obtain that the universal cover
$\widetilde{M}$
of
$M$
is an absolute Lipschitz retract and it admits a quasisymmetric
embedding into the product of 8 metric trees.
\end{abstract}
\maketitle

\section{Introduction}
A motivation of our work is a result of
Bell and Dranishnikov \cite[Theorem 1$'$]{BD}
that the asymptotic dimension of a graph-group,
whose vertex groups have finite asymptotic dimensions, is also finite.
The fundamental groups of graph-manifolds are graph-groups.
However, it is unclear whethere there exists an analogue of the
Bell-Dranishnikov result for the linearly-controlled asymptotic dimension.
We present here a construction of suitable coverings for the universal cover
of a graph-manifold that allows us to
prove the following:

\begin{thm}\label{t} For the fundamental group of a graph-manifold
$M$
taken with any word metric and its linearly-controlled
asymptotic dimension, we have
$\lasdim \pi_1(M)\leq 7$.
\end{thm}

\begin{sle}\label{sle1} Let
$\widetilde{M}$
be the universal cover of a graph-manifold
$M$.
Then
$3\leq \andim \widetilde{M}\leq 7$,
where
$\andim$
is the Assouad-Nagata dimension,
$M$
is taken with any Riemannian metric and
$\widetilde M$
with the metric lifted from
$M$.
\end{sle}

As applications of Corollary~\ref{sle1} and
\cite[Theorem~1.3, Theorem~1.5]{LS}, we obtain.

\begin{sle}\label{sle:a} Let
$\widetilde M$
be the universal cover of a graph-manifold
$M$
taken with a Riemannian distance
$d$
lifted from
$M$.
Then for all sufficiently small
$p\in(0,1)$,
there exists a bi-Lipschitz embedding of
$(\widetilde M,d^p)$
into the product of 8 metric trees, in particular
$(\widetilde M,d)$
admits a quasisymmetric embedding into the product of
8 metric trees.
\end{sle}

\begin{sle}\label{sle2} The universal cover
$\widetilde{M}$
of a graph-manifold is an absolute Lipschitz retract,
that is, given a metric space
$X$,
there is
$C>0$
such that for every subset
$A\subset X$
and every
$\lambda$-Lipschitz
map
$f\colon A\to \widetilde{M}$, $\lambda>0$,
there exists a
$C\lambda$-Lipschitz
extension $\bar{f}\colon X\to\widetilde{M}$
of
$f$.
\end{sle}

\section{Preliminaries}\label{konstr}
In this section we recall some basic definitions
and notations.
Let $X$ be a metric space.
We denote by $|xy|$
the distance between $x,y\in X$, and
$d(U,V):=\inf\{|uv|\mid u\in U,v\in V\}$
is the distance between $U,V\subset X$.
We write
$B_r(x)=\{x'\in X\mid |xx'|\le r\}$
for the ball with the center $x$ and the radius $r$.
A map $f\colon X\to Y$ is said to be
\textit{quasi-isometric} if there exist
$\lambda\geq 1,C\geq 0$ such that
$$\frac{1}{\lambda}|xy|-C\leq |f(x)f(y)|\leq \lambda |xy|+ C$$
for each $x,y\in X$.
Metric spaces $X$ and $Y$ are called
\textit{quasi-isomeric} if there is a quasi-isometric
map $f\colon X\to Y$ such that $f(X)$ is a net in $Y$.
In this case $f$ is called a \textit{quasi-isometry}.

We say that a family $\EuScript U$ of subsets of
$X$ is a \textit{covering} if for each
point $x\in X$ there is a subset $U\in \EuScript U$
such that
$x\in U$.
A family $\EuScript U$ of sets
is \textit{disjoint} if each two sets
$U,V\in \Uu$ are disjoint. The union
$\EuScript U=\cup \{\EuScript
U^{\alpha}\mid \alpha\in \EuScript A\}$
of disjoint families
$U^{\alpha}$
is said to be
$n$-{\em colored},
where
$n=|\EuScript A|$
is the cardinality of
$\EuScript A$.

Also recall that a family $\EuScript U$ is
\textit{$D$-bounded}, if the diameter of every
$U\in \EuScript U$ does not exceed $D$, $\diam U\leq D$.
An $n$-colored family of sets $\EuScript U$ is
\textit{$r$-disjoint}, if for every color
$\alpha\in \EuScript A$ and each two sets
$U,V\in \EuScript U^{\alpha}$ we have $d(U,V)\geq r$.

The linearly-controlled asymptotic
dimension is a version of the Gromov's
asymptotic dimension, $\asdim$.

\begin{dfn}(Roe \cite{Roe})
The \textit{linearly-controlled asymptotic dimension}
of a metric space $X$, $\lasdim X$, is the least
integer number $n$ such that for each sufficiently large
real $R$ there exists an $(n+1)$-colored,
$R$-disjoint, $CR$-bounded covering of the space $X$,
where the number $C>0$ is independent of $R$.
\end{dfn}

A \textit{tripod} in a geodesic metric
space $X$ is a union of three geodesic segments
$xt\cup yt\cup zt$ which have only one common
point $t$.
We say that $t$ is the \textit{center} of the tripod.
A geodesic metric space $X$ is called a
\textit{metric tree} if each triangle in it is a
tripod (possibly degenerate).

Let $T$ be a metric tree. Recall the well-known
construction of 2-colored, $r$-disjoint, $3r$-bounded
covering of
$T$
(see \cite{Roe}).
Fix some vertex $v\in T$. Suppose $k\geq 1$ is an integer
and consider the annuli
$$A_k:=\{x\in T\mid kr\leq |vx|\le (k+1)r\}$$
of the width $r$ and with the center $v$.
For $x,y\in A_k$ we say that $x\sim_k y$ if
$|zv|\ge(k-\frac{1}{2})r$, where $z$ is the center of
the tripod with vertices $x,y,v$.
One can easily check (see  \cite{Roe}) that $\sim_k$
is an equivalence relation on $A_k$.
We will make use of the following 2-colored covering:
$\EuScript U:=\EuScript U^1\cup \EuScript U^2$ where
$\EuScript U^1$ consists of the ball $B_r(v)$ and
the equivalence classes for $\sim_k$ corresponding
to all even $k$.
The family $\EuScript U^2$ consists of the equivalence
classes for $\sim_k$ corresponding to all odd $k$.
It is shown in \cite{Roe} that this covering is
$r$-disjoint and $3r$-bounded.
We call it the \textit{standard $(r,v)$-covering} of~$T$.

We also use a modified covering where the ball
$B_r(v)$ is replaced by the ball $B_d(v)$ with $0<d\leq r$.
The annuli $A_k$ are replaced by the annuli
$$B_k:=\{x\in T\mid (k-1)r+d\leq |vx|\le kr+d\}$$
of the width $r$ and with the center $v$.
Now for each integer $k>1$ as well as for
$k=1$
and
$d\ge r/2$,
we introduce on
$B_k$
the equivalence relation
$x\sim_k y$
iff
$|zv|\ge(k-\frac{3}{2})r+d$.
If
$k=1$
and
$d<r/2$,
then we say that all points in $B_1$ are
$\sim_1$-equivalent.
As before one can easily check that the modified
covering is $r$-disjoint and $3r$-bounded.
We call this covering the
\textit{standard $(r,d,v)$-covering} of $T$.\label{konstr1}

Let $G$ be a finitely generated group and $S\subset G$
a finite symmetric generating set for $G$ ($S^{-1}=S$).
Recall that a \textit{word metric} on the group $G$
(with respect to $S$) is the metric defined by the norm
$\|\cdot\|_S$ where for each $g\in G$ its norm $\|g\|_S$
is the smallest number of elements of $S$ whose product
is $g$.
In this paper we will consider only finitely generated
groups with a word metric.

\begin{dfn}
By a \textit{graph-manifold} we mean a closed three
dimensional manifold which admits a finite cover by
a connected orientable manifold that is glued from
blocks along boundary tori in a regular way.
By a \textit{block} we mean a trivial $S^1$-fibration
over a surface with negative Euler characteristic and
nonempty boundary. A gluing of blocks $M_1$ and $M_2$
(possibly $M_1=M_2$) along boundary tori
$T_1\subset\partial M_1$ and  $T_2\subset\partial M_2$
is said to be \textit{regular}, if $S^1$-fibers, coming
up from $T_1$ and $T_2$, are not homotopic on the gluing
torus.
\end{dfn}

Fix a Riemannian metric on a graph-manifold $M$.
The following fact  is well known \cite{BBI}.

\begin{lmm}\label{lem:qmetric}
Let $Y$ be a compact metric space with length metric.
Let $X$ be a universal cover of $Y$ with the metric
lifted from $Y$. Then $X$ is quasi-isometric to the
group $\pi_1(Y)$ with any word metric.
\end{lmm}

\section{Graph-manifold model}
\subsection{Description of the model}
\label{subsect:model}

Consider the hyperbolic plane
$\mathbb{H}^2_{\kappa}$
having a curvature
$-\kappa$ $(\kappa>0)$
such that the side of a
rectangular equilateral hexagon $\theta$ in the
plane $\mathbb{H}^2_{\kappa}$ has the length 1.
Let $\rho$ be the distance between the
middle points of sides, which have a common adjacent
side, $\delta$ the diameter of $\te$. It is clear
that $\rho\leq 2$, $\delta\leq 3$.

We mark each second side of $\theta$ (so we have
marked three sides) and consider a set
$\Hh$
defined as follows. Take the subgroup
$G_\theta$
of the isometry group of
$\mathbb{H}^2_{\kappa}$
generated by reflections in (three) marked sides of
$\theta$
and let
$\Hh$
be the orbit of
$\theta$
with respect to
$G_\theta$.
Then
$\Hh$
is a convex subset in
$\mathbb{H}^2_{\kappa}$ divided into hexagons that are
isometric to $\theta$.
Furthermore, the boundary of $\Hh$ has infinitely
many connected components each of which is a geodesic line in
$\mathbb{H}^2_{\kappa}$.

The graph
$T_0$
dual to the decomposition of
$\Hh$
into hexagons is the standard binary tree whose vertices
all have degree three.
Any metric space isometric to
$\Hh$
will be called a
$\theta$-{\em tree}.
Given a vertex
$p$
of $T_0$,
we denote by
$\theta_p$
the respective hexagon in
$\Hh$.

Every boundary line
$\ell$
of
$\Hh$
is the union of segments of length one each of which is
a side of a hexagon in
$\Hh$.
Thus
$\ell$
determines a sequence of vertices of
$T_0$
consecutively connected by edges which form a line
$\bar{\ell}$\label{ell}
in
$T_0$.

Consider the line $\R$ and its partition into segments
$\R=\bigcup \{[i,i+1]\mid i\in \Z\}$.
Then the graph
$S_0$\label{s0}
dual to this partition is an infinite tree with all
the vertices having degree 2, in particular,
$S_0$
is homeomorphic to
$\R$.
Given a vertex
$q$
of
$S_0$,
we denote by
$I_q\subset\Hh$
the respective segment.

\begin{zam} In what follows we will consider $T_0$ and
$S_0$ as metric spaces with metrics such that the length
of each edge is
equal to 1. Then these metric spaces are
metric trees. We will denote the set of vertices in
$T_0$ by $\V(T_0)$,
and the set of vertices in $S_0$ by $\V(S_0)$.
\end{zam}

\begin{zam} Given distinct vertices
$p,q\in \V(T_0)$
in
$T_0$
and
$x\in \te_p,y\in\te_q$,
we have
$$|pq|-1\leq |xy|\leq (|pq|-1)\rho+2\delta. \eqno(\star)$$

Similarly, given distinct
$p,q\in \V(S_0)$
and
$x\in I_p,y\in I_q$,
we have
$$|pq|-1\leq |xy|\leq |pq|+1. \eqno(\star\star)$$
\end{zam}

\begin{zam}\label{vyp} Note that for each connected
subset $U\subset T_0$ the set $U':=\cup\{\te_p\mid p\in
U\cap \V(T_0)\}$ is
convex. Similarly, for each connected subset $V\subset
S_0$ the set $V':=\cup\{I_q\mid q\in V\cap \V(S_0)\}$ is
convex.
\end{zam}

Now, we describe the following metric space
$X$.
\begin{constr}
Let $T$ be an infinite tree all whose vertices have
countably infinite degree. Denote by $\V(T)$ the set of
all the vertices of $T$ and by $\E(T)$ the set of all the edges
of $T$. For each $v\in \V(T)$ consider the metric
product $X_v:=\Hh\times \R$ and fix a bijection $b_v$
between the set $\E_v$ of all the edges adjacent to $v$
and all the boundary lines in the factor $\Hh$.
Note that the partition of $\Hh$ into hexagons and the
partition of $\R$ into segments give rise to a partition of
$X_v$ into sets isometric to $\te\times [0,1]$ which we call
\textit{bricks}.
Furthermore, each boundary plane
$\sigma\subset X_v$
is obtained as the product of a boundary
line $\ell$ in the $\te$-tree factor and the factor
$\R$,
and the partition of $X_v$ into bricks gives rise to a
partition of $\sigma$ into squares with
side length 1, which we call \textit{the grid} on the
plane $\sigma$.
Next, consider for each edge $uv$ in $\E(T)$ the
boundary planes
$\sigma_u=b_u(uv)\times\R_u\subset X_u$ and
$\sigma_v=b_v(uv)\times\R_v\subset X_v$, respectively,
and glue them by an isometry that glues the grids on $\sigma_u$ and
$\sigma_v$
while flips the factors
$\R_u$, $\R_v$
to the factors
$b_v(uv)$, $b_u(uv)$
respectively.
We call the obtained space the \textit{model} and
denote it by $X$.
\end{constr}

Recall that a complete geodesic metric space $Y$ is
called a \textit{$\cat$-space}, if the following is
satisfied:
for each triple of points $x,y,z\in X$ and a point
$t\in yz$
consider points
$\bar{x},\bar{y},\bar{z}$ on the plane $\R^2$  such
that $|xy|=|\bar{x}\bar{y}|,\ |yz|=|\bar{y}\bar{z}|,\
|xz|=|\bar{x}\bar{z}|$
and the point $\bar{t}\in [\bar{y}\bar{z}]$, such that
$|yt|=|\bar{y}\bar{t}|$, then we have
$|xt|\leq|\bar{x}\bar{t}|$.
A $\cat$-space is also called an \textit{Hadamard space}.
The spaces $\Hh$ and $\R$ are $\cat$-spaces,
so
$X_v$, $v\in\V(T)$
and $X$ are $\cat$-spaces too
(see for example \cite{BBI}).

\subsection{The model $X$ as the universal cover of a
closed graph-manifold}
Here we construct a closed graph-manifold $Q$
(with $\mathrm{C}^1$-smooth Riemannian nonpositively
curved metric) whose universal cover $\widetilde{Q}$
is isometric to $X$.

Let $\omega$ be the union of adjacent hexagons
$\te_1,\te_2\subset \Hh$.
Denote by $\tau\colon \te_1\to \te_2$ the reflection
with respect to the common side $\te_1\cap\te_2$ and
consider the surface $P_0$, obtained from $\omega$ by
gluing the marked sides matched by $\tau$.
Then $P_0$ is a surface of constant curvature $-\kappa$
with geodesic boundary.
The boundary $\partial P_0$ consist of 3 components
$\gamma_1,\gamma_2$ and $\gamma_3$ of length 2.
The hexagon vertices subdivide every $\gamma_{i},\
i=1,2,3$,
into a pair of segments of length 1.

Let us mark the images of hexagons' vertices
(so, there will be 6 marked points).
Note that the universal cover of $P_0$ is isometric to
$\Hh$.
Consider the circle $S^1_2$ of length 2 and mark two
diametrally opposite points $s_1$ and $s_2$ on it.
For each boundary component $\gamma_i$ $(i\in\{1,2,3\})$
consider an isometry $\Phi_i\colon \gamma_i\to S^1_2$
taking the marked points into the marked points.

Let $P:=P_0\times S^1_2$.
Consider an isometric copy
$P'$
of
$P$.
All the elements of $P'$ will be denoted by the same
letters but with prime. Let $Q$ be the manifold
obtained by gluing $P$ and $P'$ along the boundary tori
$\gamma_i\times S^1_2\subset \partial P$ and
$S^1_2\times \gamma'_i\subset \partial P'$
according to isometries
$$\Phi_i\times\Phi'^{-1}_i\colon\  \gamma_i\times
S^1_2 \longrightarrow S^1_2\times \gamma'_i,\ i\in\{1,2,3\}.$$

Note that
$Q$
is a graph-manifold the universal cover
$\widetilde{Q}$
of which consist of blocks isometric
to $\Hh\times\R$.
Any two of these blocks either are disjoint
or intersect over a plane covering a torus
$\gamma_i\times S^1_2$, $i\in\{1,2,3\}$.
This plane is divided by the lines covering the
circles $\{x_i\}\times S^1_2$ and
$\gamma_i\times \{s_i\}$ (where $x_i\in\gamma_i$ and
$s_i\in S^1_2$ are marked points) into squares with
side 1.
It is clear that $\widetilde{Q}$ isometric to $X$.
We recall the following result
\cite[Theorem 2.1]{BN}\label{bilip}.

\begin{thm}\label{teorem}
Any two graph-manifolds have bi-Lipschitz homeomorphic
universal covers, in particular, their fundamental
groups
are quasi-isometric.
\end{thm}

It follows from this theorem that for any
graph-manifold $M$ its fundamental group $\pi_1(M)$
is quasi-isometric to the fundamental group $\pi_1(Q)$
of the manifold $Q$ which is by Lemma~\ref{lem:qmetric}
quasi-isometric to the model $X$.
Therefore, since $\lasdim$ is a quasi-isometry invariant
(see e.g. \cite{LS}), we have $\lasdim \pi_1(M)=\lasdim X$.
So, to prove  Theorem \ref{teorem} it suffices to show that
$\lasdim X\leq 7$.
To this end, we construct for each $R\in\N$, $R>10$
an $R$-disjoint, $CR$-bounded 8-colored covering of the
model $X$, where $C=88$.

\section{Coverings of the model}
\label{sect:cover_model}
We briefly describe the construction of the required
covering of
$X$.
We fix
$R\in\N$, $R>10$.
First, we construct for every block
$X_u$, $u\in \V(T)$,
a 4-colored,
$2R$-disjoint,
$15R$-bounded
covering
$\Ww_u$
which is the product of standard 2-colored coverings
of the factors in the decomposition
$X_u=\Hh\times\R$.
For different neighboring blocks these covering are
compatible (see Lemma~\ref{zz}) and form together 4-colored covering
$\widetilde\Ww$
of
$X$.
The covering
$\widetilde\Ww$
is
$R$-disjoint
(see Proposition~\ref{u1}), and this is the main result
of the current section. However, all the members of
$\widetilde\Ww$
have infinite diameter. Fortunately, every member
$W$
of
$\widetilde\Ww$
has a tree-like structure, and we construct in
sect.~\ref{susect:adcolors} a 2-colored covering of
$W$
similar to the standard 2-colored covering of a tree.
Typical
$W$
is disconnected having every connected component convex.
For technical reasons, it is convenient first to extend every
$W\in\widetilde\Ww$
to a convex
$\overline W\supset W$
and then to cover
$\overline W$
by a 2-colored
$\Xx_W$.
The extension procedure is described in sect.~\ref{subsect:extension}.
The convexity plays an important role in proving the properties of
$\Xx_W$.
After that, we forget about the extension
$\overline W$
and restrict
$\Xx_W$
to
$W$.
In that way, we produce the required 8-colored
covering
$\Xx$
of
$X$.

\subsection{Construction of a covering of a
block}\label{sss}

Put
$N:=2R+1$.
Consider an arbitrary block $X_v=\Hh\times \R\subset X$.
We choose a brick $\te_p\times I_q\subset X_v$,
where $p\in\V(T_0)$ and $q\in\V(S_0)$.
Consider the standard $(N,p)$-covering
$\Uu=\Uu^1\cup\Uu^2$
of the tree $T_0$~(see~section~\ref{konstr}).
We build a 2-colored covering
$\widetilde{\Uu}=\widetilde{\Uu}^1\cup\widetilde{\Uu}^2$
of the $\te$-tree as follows:
for each set $U\in\Uu^i$, $i\in\{1,2\}$ we define the set
$$\widetilde{U}:=\bigcup\{\te_{p'}\mid p'\in U\cap\V(T_0)\}$$
and consider the families $$\widetilde{\Uu}^i:=\{\widetilde{U}\mid U\in \Uu^i\},\ i\in \{1,2\}.$$
It is clear that the family $\widetilde{\Uu}$ is a
2-colored covering of the
$\te$-tree. Recall that $\rho$ is the distance between
the middle points of
two sides of the hexagon $\te$, which have a common
adjacent side, and $\delta$ is its diameter.

Since the covering $\Uu$ is $N$-disjoint
and $3N$-bounded,
the inequality $(\star)$ implies that the covering
$\widetilde{\Uu}$ is $2R$-disjoint and 
$((6R+2)\rho+2\delta)$-bounded. Since $\rho\leq 2,\delta\leq 3$ and $R>10$
we have $(6R+2)\rho+2\delta\leq13R$ so
$\widetilde{\Uu}$ is $2R$-disjoint and $13R$-bounded.

Fix a natural $d\leq N$. Consider the standard
$(N,d,q)$-covering $\Vv=\Vv^1\cup\Vv^2$ of
the tree
$S_0$ (see section \ref{konstr1}).
We build a 2-colored covering
$\widetilde{\Vv}=\widetilde{\Vv}^1\cup\widetilde{\Vv}^2$
of the line $\R$ as follows: for each set
$V\in\Vv^i$, $i\in\{1,2\}$ consider the set
$$\widetilde{V}:=\bigcup\{I_{p'}\mid p'\in V\cap \V(S_0)\}$$
and the families $$\widetilde{\Vv}^i:=\{\widetilde{V}\mid V\in \Vv^i\},\ i\in \{1,2\}.$$
It is clear that the family $\widetilde{\Vv}$ is a
2-colored covering of $\R$.
Note, that since the family $\Vv$ is $N$-disjoint and $3N$-bounded the inequality
$(\star\star)$ implies that the covering
$\widetilde{\Vv}$ is
$2R$-disjoint and $(6R+2)$-bounded. Since $R>10$ we have
$6R+2\leq 7R$, so the covering $\widetilde{\Vv}$ is
$2R$-disjoint and $7R$-bounded.

Consider the 4-colored covering
$$\Ww:=\bigcup\{\Ww^{(i,j)}\mid (i,j)\in \{1,2\}\times \{1,2\}\}$$
of the block $X_v$ where
$$\Ww^{(i,j)}:=\{U\times V\mid U\in \widetilde{\Uu}^i, V\in \widetilde{\Vv}^j\}.$$
Since the coverings $\widetilde{\Uu}$ and
$\widetilde{\Vv}$ are $2R$-disjoint the covering $\Ww$
is also $2R$-disjoint, and since the covering
$\widetilde{\Uu}$ is $13R$-bounded and $\widetilde{\Vv}$
is $7R$-bounded, the covering $\Ww$ is
$R\sqrt{7^2+13^2}$-bounded, thus it is $15R$-bounded.
We will call this covering of the block $X_v$ the
\textit{$(N,d)$-covering with
the initial brick $\te_p\times I_q$}.

Consider the block
$X_v$
with the fixed initial brick
$\te_{p_v}\times I_{q_v},\ p_v\in \V(T_0),q_v\in\V(S_0)$,
and its boundary plane
$\sigma=\ell\times\R_v$,
where
$\ell$
is a boundary line of the
$\te$-tree
factor of
$X_v$.
Let $\bar{\ell}$ be the line in the tree
$T_0$
that corresponds to
$\ell$
(see section \ref{ell}),
$p'\in\bar{\ell}$
the nearest to
$p_v$
vertex. We denote
$n:=|p_vp'|$
and put
$k\equiv n~(\mod~N),~0\leq k< N$.

Let
$X_u$
be the block glued to
$\sigma$, $1\le d\le N$
a natural number,
$I_{q_u}\times \te_{p_u}$
the brick in
$X_u$
having the common boundary square with the brick
$\te_{p'}\times I_{q_v}$.
We represent the boundary plane
$\sigma$ as the product of the line $\R$
and an appropriate boundary line $\ell'$ in the
$\te$-tree factor of the block $X_u$,
$\sigma=\R_u\times\ell'$.
Consider the line
$\bar{\ell'}$ in the tree $T_0$ which corresponds to the
line $\ell'$.
Then the vertex $p_u$ belongs to
$\bar{\ell'}$.
We choose a vertex $p'_u$ of $T_0$ such that the vertex
$p_u$ is the nearest to $p'_u$ point of the line
$\bar{\ell'}$ and the distance $|p_up'_u|$ is equal to
$N-d$.

\begin{lmm}\label{zz}
The $(N,d)$-covering $\Ww_v$ with the initial brick
$\te_{p_v}\times I_{q_v}$ of the block $X_v$ and the
$(N,N-k)$-covering $\Ww_u$ with the initial brick
$I_{q_u}\times \te_{p'_u}$ of the block $X_u$
\textit{agree} on the common boundary plane $\sigma$ of
the blocks $X_u$ and $X_v$, i.e. for each pair of numbers
$(i,j)\in \{1,2\}\times \{1,2\}$ any two sets
$V\in \Ww^{(i,j)}_v$ and $U\in \Ww^{(j,i)}_u$ either
are disjoint or their intersections with $\sigma$
coincide.
\end{lmm}

\begin{proof}
Indeed since the coverings $\Ww_v$ and $\Ww_u$ are
products of the coverings of the corresponding factors, it
is enough to check that the coverings of
the factors agree. We have
$\sigma=\ell\times\ell'$
and the
$\R$-factor of the block $X_u$ is glued to the
line $\ell$.
Let $\bar{\ell}$ be the line corresponding to the line
$\ell$ in the tree $T_0$.
If we consider the line $\bar{\ell}$ as an isometric
copy of the tree $S_0$ then the standard
$(N,p_v)$-covering of the tree $T_0$
induces a covering $\Uu'$ on the line $\bar{\ell}$
which is a standard $(N,N-k,p')$-covering. Therefore
since the side of the hexagon $\te_{p'}$ lying on
$\ell$ is glued to the segment $I_{q_u}$ of the
$\R$-factor of the block
$X_u$, the standard $(N,N-k,q_u)$-covering of this
$\R$-factor coincides with the covering $\Uu'$.
The case of the other factors is similar.
\end{proof}

\subsection{Constructing compatible coverings of blocks}
Recall that
$T$
is the canonical simplicial tree with all the vertices
having the infinite countable degree. We consider
the metric on
$T$
with respect to which all the edges have length one.
We fix an arbitrary vertex $v\in \V(T)$ and
call it the \textit{root}
of the tree
$T$.
We define the {\em level} function
$l:\V(T)\to\Z$
by
$l(v'):=|vv'|$.
We build coverings of blocks of the type
discussed earlier which agree with the
coverings of the neighbor blocks on the common boundary
planes (see sect.~\ref{sss})
by induction on the level $m$.

The base of induction : $m=0$.\\
Let $\te_p\times I_q$ be an arbitrary brick in $X_v$.
Consider the $(N,N)$-covering with the initial brick
$\te_p\times I_q$ of the block $X_v$.
We denote this covering by $\Ww_v$.

The induction step: $m \rightarrow m+1$.\\
Suppose that we have already built the coverings
$\Ww_{v'}$ which are compatible for all the blocks
$X_{v'}$ such that $l(v')\leq m$. Consider
an arbitrary vertex $u\in \V(T)$ such that $l(u)=m+1$.
The block $X_u$ has exactly one boundary plane  $\sigma$
along which it is glued with a block
$X_{w}$ such that $l(w)=m$.
$X_{w}$ is already covered by $\Ww_{w}$ which is an
$(N,d_w)$-covering with the initial brick
$\te_{p_w}\times I_{q_w}$. Consider the plane
$\sigma=\ell\times\R$ as the product of a boundary line
$\ell$ in the $\te$-tree and the factor $\R$.
Let $\bar{\ell}$ be the line in
$T_0$ corresponding to $\ell$ (see section \ref{ell}).
We take the vertex $p'$ in
$\bar{\ell}$ which is the nearest to $p_w$
and denote
$n_w:=|p_wp'|$. Let $d_u\equiv n_w~(\mod~ N),~0\leq d_u<N$.
Let $I_{q_u}\times \te_{p_u}$ be a brick in  $X_u$
having a common boundary square
with the brick $\te_{p'}\times I_{q_w}$.
We represent the boundary plane $\sigma$ as the product of
the line $\R$ and a boundary line in the $\te$-tree
factor of the block $X_u$, $\sigma=\R\times\ell'$.
Consider the line $\bar{\ell'}$ in $T_0$ which
corresponds to the line $\ell'$.
It is clear that the vertex $p_u$ belongs to
$\bar{\ell'}$.
We choose a vertex $p'_u$ in $T_0$ such that the vertex
$p_u$ is the nearest to $p'_u$ point on the line
$\bar{\ell'}$ and the distance $|p_up'_u|$ is equal to
$N-d_w$.
We denote by $\Ww_u$ the $(N,N-d_u)$-covering of the block
$X_u$  with the initial brick $I_{q_u}\times \te_{p'_u}$.
It follows from Lemma \ref{zz} that the coverings
$\Ww_w$
and
$\Ww_u$
of the blocks
$X_w$
and
$X_u$
respectively agree on the boundary plane $\sigma$.
$\Ww_u$ is the covering of the type we need so the
induction step is completed.

\subsection{Constructing a covering of the model}
\label{subsect:cover_model}
Consider the following 4 families of subsets of the
model $X$. For each pair of numbers
$(i,j)\in\{1,2\}\times\{1,2\}$ we denote
$$\Ww^{(i,j)}:=\{U\mid U\in \Ww^{(i,j)}_{u},\ l(u)\hbox{ is even}\}
\cup \{U\mid U\in \Ww^{(j,i)}_{u},\ l(u)\hbox{ is odd}\}.$$

Moreover, for each pair of numbers
$(i,j)\in\{1,2\}\times\{1,2\}$ we consider
the following  relation $\sim_{(i,j)}$ on the family
$\Ww^{(i,j)}$: $U\sim_{(i,j)} V$ iff there are sets
$U_0=U,U_1,\ldots,U_{n+1}=V$ in the family $\Ww^{(i,j)}$
such that
$U_j\cap U_{j+1}\ne\emptyset$
for all
$j\in \{0,\ldots,n\}$.
It is clear that
$\sim_{(i,j)}$
is an equivalence relation.

\begin{utv}\label{u1}
Suppose $W,U\in \Ww^{(i,j)}$ and $d(W,U)<R$. Then
$W\sim_{(i,j)} U$.
\end{utv}
We need two facts (see Lemma~\ref{l1} and \ref{l2})
for the proof of Proposition~\ref{u1}.

\begin{lmm}\label{l1} Given
$U\in \Ww^{(i,j)}$, $U\subset X_u$,
assume that a boundary plane
$\sigma\subset X_u$
meets
$U$, $U\cap\sigma\ne\emptyset$,
and separates a point
$x\in X$
and
$U$.
Then
$d(x,U)=d(x,U\cap\sigma)$.
\end{lmm}

\begin{proof} Without loss of generality, we assume that
$x\not\in\sigma$.
Then, since the hyperplane
$\sigma\subset X$
is convex, the segment
$xy$
meets
$\sigma$
over a point for every
$y\in U$, $xy\cap\sigma=t$.

Recall that we represent the block
$X_u$
as the product
$X_u=\Hh\times\R$
and also
$U=\widetilde U_\te\times\widetilde U_\R$,
where
$\widetilde U_ \te $, $\widetilde U_\R$
are members of the appropriate coverings of
$\Hh$
and
$\R$
respectively. Then
$t=(t_{\te},t_{R})\in\ell\times\R=\sigma$,
where
$\ell\in\Hh$
is a boundary line, and
$y=(y_\te,y_R)$,
where
$y_\te\in\widetilde U_\te$, $y_R\in\widetilde U_\R$.
For
$t'=(t_\te,y_R)\in\sigma$,
we have
$$|xt'|\le|xt|+|tt'|\le|xt|+|ty|=|xy|.$$
If
$t_\te\in\widetilde U_\te$,
then
$t'\in U\cap\sigma$
because
$y_R\in\widetilde U_R$,
and hence
$|xy|\ge d(x,U\cap\sigma)$.
Thus in what follows we assume that
$t_\te\not\in\widetilde U_\te$.

Let $\te_p\times I_q$
be the initial brick in
$\Ww^{(i,j)}_u$.
Passing to the tree
$T_0$,
we denote by $U_0$ the set of all vertices $q'$ in
$\V(T_0)$ such that
$\te_{q'}\subset\widetilde{U}_{\te}$.
Then by the construction of
$\Ww^{(i,j)}_u$, $U_0$
is the vertex set of a member of the standard
$(N,p)$-covering of the tree $T_0$.
Denote this member by
$U'$.
There is
$k\in\N\cup \{0\}$
such that $U'$ lies in the
$N$-annulus
$$A_k=\{q\in T_0\mid kN\leq|qp|\le(k+1)N\}$$
centered at
$p$.
We denote by
$\bar\ell$
the line in
$T_0$
corresponding to
$\ell$
and pick a vertex
$p_0\in\bar\ell$
of
$T_0$
such that
$t_\te\in\te_{p_0}$.
There are at most two such vertices, and we assume
that
$p_0\in U'$
if
$t_\te\in\widetilde U_\te$.
Then the assumption
$t_\te\not\in\widetilde U_\te$
is equivalent to
$p_0\not\in U'$.

Since
$y_\te\in\widetilde U_\te$,
there is a vertex
$p_y$
of
$T_0$
such that
$y_\te\in\te_{p_y}$
and
$p_y\in U'$.
We show that
$\dist(p_0,\bar\ell\cap U')\le|p_0p_y|$
in
$T_0$.
The (unique) point
$p'\in\bar\ell$
closest to
$p$
separates
$\bar\ell$
into the rays
$\bar\ell_i$, $i=1,2$.
Each of the rays
$\bar\ell_1$, $\bar\ell_2$
intersects the annulus
$A_k$
over a segment, and at least one of
them lies in
$U'$.
We assume without loss of generality that
$\bar\ell_1\cap A_k\subset U'$.
There is a shortest segment
$\bar\gamma\subset\bar\ell$
between
$p_0$
and
$\bar\ell\cap U'$.
If
$p_0\in\bar\ell_1$,
then
$\bar\gamma\subset\bar\ell_1$
is a shortest path in
$T_0$
between
$p_0$
and
$U'$.
Thus the length
$|\bar\gamma|=d(p_0,U')$
and hence
$d(p_0,\bar\ell\cap U')\le|p_0p_y|$.

It remains to consider the case
$p_0\in\bar\ell_2$
and
$\bar\ell_2$
is disjoint with
$U'$
(the last happens exactly when
$k>1$
and
$|pp'|\le(k-\frac{1}{2})N)$.
As above, we have
$|p'p_y|\ge d(p',U')=d(p',\bar\ell\cap U')$.
On the other hand, any path in
$T_0$
between
$p_0$
and
$U'$
passes over
$p'$.
Thus
$|p_0p_y|=|p_0p'|+|p'p_y|\ge|\bar\gamma|=d(p_0,\bar\ell\cap U')$.

The geodesic segment
$ty\subset X_u$
projects to the geodesic segment
$t_\te y_\te$
in the factor
$\Hh$
of
$X_u$.
While moving from
$t_\te$
to
$y_\te$
along
$t_\te y_\te$,
one passes the segment
$p_0p_y\subset T_0$
and
$|t_\te y_\te|\ge|p_0p_y|-1$
by
$(\star)$.

Let
$\gamma=t'z\subset\ell\times y_R\subset\sigma$
be the segment between
$t'$
and
$z\in U\cap\sigma$
that corresponds to
$\bar\gamma\subset\bar\ell$.
The length of
$\gamma$
is at most the number of {\em interior} vertices of
$\bar\gamma$
plus one, thus
$|\gamma|\le|\bar\gamma|-1$.
Since
$|p_0p_y|\ge|\bar\gamma|$,
we obtain
$|t_\te y_\te|\ge|\gamma|$.
It follows
$|ty|\ge|tz|$
because both the triangles
$tt'y$
and
$tt'z$
are flat with right angles at
$t'$.
Then moreover
$|xy|\ge|xz|\ge d(x,U\cap\sigma)$.
\end{proof}

\begin{lmm}\label{l2} Suppose
$\sigma$
is the common plane of the blocks
$X_u$
and
$X_{u'}$
such that
$|vu'|<|vu|$,
where
$v$
is the root of
$T$.
Assume that
$U\in \Ww^{(i,j)}$, $U\subset X_u$,
is disjoint with
$\sigma$, $U\cap\sigma=\emptyset$.
Then
$d(U,\sigma)\geq R$.
\end{lmm}

\begin{proof} Let
$\te_p\times I_q$
be the initial brick in
$\Ww^{(i,j)}_u$.
Recall that
$U=\widetilde{U}_{\te}\times\widetilde{U}_R$,
where
$\widetilde{U}_{\te}$ and $\widetilde{U}_R$
are members of the appropriate coverings of the
$\te$-tree
$\Hh$
and
$\R$
respectively. Denote by
$U_0$
the set of all vertices
$q'\in\V(T_0)$
such that the hexagon
$\te_{q'}\in \widetilde{U}_{\te}$.
Then
$U_0$
is the vertex set of a member
$U'$
of the standard
$(N,p)$-covering
of the tree $T_0$.
There is
$k\in \N\cup \{0\}$
such that $U'$
is a subset of the
$N$-annulus
$$A_k=\{q\in T_0\mid kN\leq|qp|\le(k+1)N\}$$
centered at
$p$.

Representing
$\sigma=\ell\times\R$,
where
$\ell\subset\Hh$
is an appropriate boundary line, we let
$\bar\ell\subset T_0$
be the line corresponding to
$\ell$.
It follows from the assumption on
$\sigma$
and the construction of the family
$\Ww^{(i,j)}$
that
$d(p,\bar\ell)\le N$
and thus
$\bar\ell$
meets every annulus
$A_k$, $k\ge 0$.
By the assumption,
$\ell$
is disjoint with
$\widetilde U_\te$,
hence
$\widetilde U_\te$
is disjoint with any hexagon of
$\Hh$
that meets
$\ell$.
Thus
$\bar\ell$
is disjoint with
$U'$.
Since
$\bar\ell\cap A_k\neq\emptyset$
is covered by members of the standard
$(N,p)$-covering
of
$T_0$
having the same color
$i\in\{1,2\}$
as
$U'$
has and the covering is
$N$-disjoint,
we have
$d(U',\bar\ell\cap A_k)\ge N$.
Therefore
$d(U',\bar\ell)\ge N/2$.
Moreover,
$N=2R+1$
is odd and the distance
$d(U',\bar\ell)$
is integer, thus we actually have
$d(U',\bar\ell)\ge (N+1)/2=R+1$.
Using $(\star)$,
we obtain
$d(\widetilde U_\te,\ell)\ge R$
and thus
$d(U,\sigma)\ge R$.
\end{proof}

\begin{proof}[Proof of Proposition~\ref{u1}]
Suppose that $U$ lies in the block $X_u$ and $W$ lies in
the block $X_w$.
We prove the proposition by induction on
$m=l(u)+l(w)$.

The base of induction : $m=0$.\\
It follows from the assumption
$d(U,W)<R$
that
$U=W$,
because different members of the family
$\Ww^{(i,j)}$
are
$N$-disjoint.

The induction step: $m \rightarrow m+1$.\\
We assume that the assertion is true for all
$l(u)+l(w)\le m$
and consider thå case
$l(u)+l(w)=m+1$.
Without loss of generality suppose
$l(u)\leq l(w)$.
Then
$l(w)\ge 1$.
Thus there exists a boundary plane $\sigma$ in $X_w$
such that
$l(w')=l(w)-1$
for the block
$X_{w'}$,
which is glued to
$X_w$
over
$\sigma$.

We have
$d(W,\sigma)<R$
because
$U$, $W$
are separated by
$\sigma$.
Therefore
$W\cap\sigma\ne\emptyset$,
since
$X_w$
and
$\sigma$
satisfy the conditions of Lemma~\ref{l2}.
By Lemma~\ref{l1},
$d(x,W)=d(x,W\cap\sigma)$
for every
$x\in U$,
thus
$d(U,W\cap\sigma)<R$.
There is
$W'\in\Ww^{(i,j)}$, $W'\subset X_{w'}$,
such that
$W'\cap\sigma=\sigma\cap W$
by construction of the family
$\Ww^{(i,j)}$,
in particular,
$W'$
is equivalent to
$W$.
Then
$d(W',U)<R$,
and thus
$W'$
is equivalent to
$U$
by the inductive assumption. Hence
$W$
is equivalent to
$U$.
\end{proof}

Given
$(i,j)\in\{1,2\}^2$,
for each equivalence class
$\widetilde{W}$
of the equivalence relation
$\sim_{(i,j)}$
on the family
$\Ww^{(i,j)}$
consider the set
$W:=\bigcup\{U\mid U\in \widetilde{W}\}$.
We denote by
$\widetilde{\Ww}^{(i,j)}$ the family of all such sets.
Proposition~\ref{u1} implies that the covering
$$\widetilde{\Ww}
=\bigcup\{\widetilde{\Ww}^{(i,j)}\mid (i,j)\in\{1,2\}^2\}$$
of
$X$
is 4-colored and
$R$-disjoint.
However,
$\widetilde\Ww$
is unbounded, since the diameter of each of its members is infinite.

\section{Proof of Theorem~\ref{t}}
\subsection{Extension of the covering members}
\label{subsect:extension}
Let
$W$
be a member of the covering $\widetilde{\Ww}$.
The sets
$W_u=X_u\cap W$, $u\in\V(T)$,
while nonempty, form the equivalence class
$\widetilde{W}$
which corresponds to
$W$.
We have
$W_u=U\times V$,
where
$U\subset\Hh$
and
$V\subset\R$
are members of the coverings $\widetilde{\Uu}^i$ and
$\widetilde{\Vv}^j$ respectively for some pair
$(i,j)\in \{1,2\}\times \{1,2\}$.

Recall that
$U$
is the union of the hexagons labelled by all the vertices of the tree
$T_0$
(see sect.~\ref{subsect:model}
for definition of
$T_0$ )
lying in a member
$U'$
of the standard
$(N,p)$-covering
of
$T_0$.
There is
$k\in\N\cup\{0\}$
such that
$U'$
is contained in the
$N$-annulus
$$A_k=\{q\mid kN\leq|qp|\le(k+1)N\}\subset T_0$$
centered at
$p$.

Now, we construct an extension
$\overline U'$
of
$U'$
as follows. For each
$x\in U'$
consider the geodesic segment
$\gamma_x=x'x\subset px$
of length
$N/2$
(if
$|px|<N/2$, then
$k=0$, $U'=B_N(p)$,
and we define
$\overline{U}':=U'$).
We put
$$\overline{U}':=\bigcup\{\gamma_x \mid x\in U'\}\subset T_0$$
and define
$$\overline{U}
:=\bigcup\{\te_v \mid v\in \overline{U}'\cap\V(T_0)\}\subset\Hh.$$
Note that the set $\overline{U}'$ is connected.
Indeed, for each
$x,y\in \overline{U}'$
there exist
$x_1,y_1\in U'$
such that
$x\in\gamma_{x_1},y\in \gamma_{y_1}$.
Let $s$ be the center of the tripod $x_1,y_1,p$.
Definition of the set $U'$ implies $|ps|\geq (k-\frac{1}{2})N$.
There are
$x'\in px_1$
and
$y'\in py_1$
such that
$|px'|=|py'|=kN$.
Then
$x_1x'$, $y_1y'\subset U'$
and either
$s\in U'$
or
$|sx'|$, $|sy'|\le N/2$.
In each case
$s$
is connected with
$x$
and
$y$
by paths in
$\overline U'$,
thus the set
$\overline{U}'$
is connected. Then
$\overline U$
is convex by Remark~\ref{vyp}.

Recall that
$V$
is the union of the intervals labelled by
all the vertices of the tree
$S_0$
(homeomorphic to
$\R$, see sect.~\ref{s0} for definition of $S_0$)
lying in a member
$V'\subset S_0$
of the standard
$(N,d,q)$-covering
of
$S_0$
for some natural
$d\le N$.

We similarly construct an extension
$\overline V'$
of
$V'$.
For each
$x\in V'$
consider the segment
$\gamma_x:=x'x\subset qx$
of length
$N/2$
(if $|qx|<N/2$, then either $V'=B_d(q)$ and we define
$\overline V':=V'$,
or
$V'=B_{d+N}(q)\setminus\int B_d(q)$ and $d<N/2$,
then we define
$\overline V':=B_{d+N}(q)$).
Consider the set
$$\overline V':=\bigcup\{\gamma_x \mid x\in V'\}\subset S_0$$
and define
$$\overline{V}:=\bigcup\{I_v \mid v\in \overline{V}'\cap\V(S_0)\}
\subset\R.$$
The set $\overline{V}'$
being a segment in
$S_0$
is connected. Then
$\overline V$
is convex by Remark~\ref{vyp}.

Therefore
$\overline W_u:=\overline U\times \overline V$
is a convex subset of
$X_u$.
Moreover,
$\overline{W}_u$
is compact by construction. We call it the \textit{extension}
of the set
$W_u$.
Note that similarly to Lemma~\ref{zz}
the extensions of the sets
$W_u$
and
$W_w$
agree for the blocks
$X_u$
and
$X_w$
having the common boundary plane
$\sigma$,
i.e. their intersections with
$\sigma$
coincide. We call the \textit{extension} of the set
$W$
the union of the extensions of all its block components
$W_u$.

Recall well known fact that for each point
$x$
in an Hadamard space
$X$
and for each closed convex subset
$A\subset X$
there exists the unique point
$p_x\in A$
such that
$|xp_x|=d(x,A)$
(the point
$p_x$
is called the metric projection
$x$
to the set
$A$).
The map sending
$x\in X$
to its metric projection
$p_x\in A$
is called the {\em projection} to
$A$
and is a
$1$-Lipschitz map.

\begin{utv} For each member
$W$ of the covering
$\widetilde{\Ww}$
its extension
$\overline{W}$
is convex.
\end{utv}

\begin{proof} Given
$x$, $y\in \overline{W}$,
we prove that the geodesic segment
$xy\subset \overline{W}$.
Let
$X_u$, $X_v\subset X$ ($u,v\in\V(T)$)
be the blocks containing
$x$, $y$
respectively. If
$u=v$,
then
$x$, $y\in\overline W_u$
that is convex. Thus we assume that
$u\neq v$.
The segment
$xy$
intersects consecutively the
boundary planes $\sigma_1,\ldots,\sigma_n$ of the blocks
$X_0=X_u,\ldots,X_n=X_w$ ($\sigma_i=X_{i-1}\cap X_i$ for
each $i\in\{1,\ldots, n\}$) at
$z_i\in \sigma_i,\ i\in\{1,\ldots, n\}$.
The sets
$\overline{W}_i:=X_i\cap \overline{W}$ are convex for
all $i\in\{0,\ldots, n\}$.
Note that for each $i\in\{1,\ldots, n\}$ the
projection of $z_i$ to
$\overline{W}_{i-1}$ and to $\overline{W}_i$ coincide.
Indeed, let
$t\in \overline{W}_i$
be the nearest point to
$z_i$.
By construction of
$\overline{W}_i$,
the metric projection of
$\overline{W}_i$
to
$\sigma_i$
coincides with
$\overline{W}_i\cap\sigma_i$,
in particular, the projection
$t_i$
of
$t$
to
$\sigma_i$
belongs to
$\overline{W}_i$.
Since
$\sigma_i$
is convex in the Hadamard space
$X$,
we have
$|z_it_i|\leq |z_it|$, therefore the nearest
point in $\overline{W}_i$ to $z_i$
lies in
$\overline{W}_i\cap\sigma_i$.
By similar argument with $\overline{W}_{i-1}$,
we see that projections of $z_i$
to
$\overline{W}_{i-1}$ and to $\overline{W}_i$ coincide.
Denote this point by
$z'_i$.
By convexity of
$\overline W_i$,
we have
$|z'_iz'_{i+1}|\le|z_iz_{i+1}|$
for every
$i=0,\dots,n$,
where
$z_0=z_0'=x$, $z_{n+1}=z_{n+1}'=y$,
i.e. the curve
$\gamma:=xz'_1\cup\ldots\cup z'_ny$
has length less or equal to that of
$xy$,
hence
$\gamma=xy$.
Thus
$xy\subset \overline{W}$.
\end{proof}

\subsection{Additional colors}
\label{susect:adcolors}

Consider the extension $\overline{W}$ of a member $W$ of
the covering $\widetilde{\Ww}$.
Note that for each block $X_u$ such that
$X_u\cap W\ne \emptyset$, the set
$\overline{W}_u:=X_u\cap \overline{W}$ is convex,
compact and
$17R$-bounded.
The last follows from the fact that
$W_u=X_u\cap W$
is
$15R$-bounded
and each point
$x\in \overline{W}_u$
lies at the distance at most
$R$
from $W_u$.

\begin{zam} Given
$u,w\in V(T)$,
any segment
$xy\subset X$
with
$x\in W_u$, $y\in W_w$
intersects
$W_{u'}$
iff the segment
$uw\subset T$
contains
$u'$.
\end{zam}

Put $c=18$.
Similarly to the standard coverings of a metric tree (see
section \ref{konstr}), we construct for every member
$W$
of the covering
$\widetilde\Ww$
a 2-colored,
$R$-disjoint,
$c'R$-bounded
covering of the extension $\overline{W}$, where $c':=5c-2$.

Given
$x\in X$, $u\in V(T)$,
we use notation
$x\overline W_u$
for the geodesic segment
$xx'\subset X$,
where
$x'$
is the metric projection of
$x$
to
$\overline W_u$.
This is well defined because
$\overline W_u$
is convex and thus
$x'$
is unique. We fix
$v\in\V(T)$
such that
$W\cap X_v\neq\emptyset$.
For every
$k\in \N\cup\{0\}$
consider the ``annulus''
$$A_k:=\bigcup\{\overline{W}_u
\mid ckR\leq
d(\overline{W}_v,\overline{W}_u)\le c(k+1)R\}\subset\overline W.$$
We define a relation
$\sim_k$ on the set $A_k$
saying that
$x\sim_0y$
for every
$x$, $y\in A_0$,
and for each natural $k$
we put
$x \sim_k y$
iff there exists a vertex
$u\in \V(T)$
such that the geodesic segments
$x\overline{W}_v$, $y\overline{W}_v$
intersect the set
$\overline{W}_u$
and
$$d(\overline{W}_u,\overline{W}_v)\geq (c(k-1)+1/2)R.$$
Note that
$x\sim_ky$
for each
$x$, $y\in\overline W_u\subset A_k$.

We prove that $\sim_{k}$ is an equivalence relation.
Indeed, suppose $x \sim_{k} y$ and $y \sim_{k} z$.
Given
$x\in \overline{W}_{v_{x}}, y\in \overline{W}_{v_{y}},$
$z\in \overline{W}_{v_{z}}$,
there are
$u\in \V(T)$ from the definition of
$\sim_{k}$ for the pair $\{x,y\}$
and
$w\in \V(T)$ for the pair $\{y,z\}$. Note that
$u$, $w$
lie on the segment
$v_yv\subset T$.
Without loss of generality suppose that $|uv|\leq |wv|$.
Then the segment
$z\overline{W}_{v}$
intersects the set $\overline W_u$,
so the points $x$ and $z$
are equivalent too.

Consider the following 2-colored covering of the set
$\overline{W}$: the sets of the first color are all the
equivalence classes of $\sim_k$, where $k$ is even, the
sets of the second color are all the equivalence classes of
$\sim_k$, where $k$ is odd.

We prove that this covering of
$\overline W$
is
$c'R$-bounded.
Indeed, assume that
$x\in \overline{W}_{v_{x}}\subset A_k$, $y\in \overline{W}_{v_{y}}\subset A_k$
are equivalent, and let
$u\in \V(T)$
be a vertex such that the segments
$x\overline W_v$, $y\overline W_v$
intersect the set
$\overline{W}_u$
and
$$d(\overline{W}_u,\overline{W}_v)\geq (c(k-1)+1/2)R.$$
For
$z$, $t\in\overline W_u$
with
$z\in x\overline W_v$, $t\in y\overline W_v$
we have
$|z\overline W_v|$, $|t\overline W_v|\ge (c(k-1)+1/2)R$.
Then
$|xz|=|x\overline W_v|-|z\overline W_v|\leq (2c-\frac{1}{2})R$,
and similarly
$|yt|\leq (2c-\frac{1}{2})R$.
Since
$\overline{W}_u$
is
$(c-1)R$-bounded,
we obtain
$$|xy|\leq |xz|+|zt|+|ty|\leq (5c-2)R=c'R.$$

Now, we prove that the constructed covering of
$\overline W$
is
$R$-disjoint.
Let $U_1\subset A_k,U_2\subset A_l$ be sets of the same
color.

Assume first that
$|k-l|\geq 2$.
Without loss of generality suppose
$k\geq l+2$.
If
$d(U_1,U_2) < R$,
then
$|xy|<R$
for some
$x\in\overline{W}_{v_x}\subset U_1$, $y\in\overline{W}_{v_y}\subset U_2$.
Note that
$d(\overline{W}_{v},\overline{W}_{v_{y}})\le c(l+1)R$,
thus $|zt|\le c(l+1)R$
for some
$z\in \overline{W}_{v_{y}}$, $t\in\overline{W}_{v}$.
Then
$$|xt|\leq |xy|+|yz|+|zt|<R+(c-1)R + c(l+1)R=c(l+2)R\le ckR.$$
This contradicts
$ckR \leq d(\overline{W}_{v},\overline{W}_{v_{x}})$.

Assume now that
$k=l$, $\overline W_a\subset U_1$, $\overline W_b\subset U_2$.
It suffices to prove that
$d(\overline W_a,\overline W_b)\ge R$.
Let
$u\in\V(T)$
be the center of the tripod
$abv\subset T$.
Then any geodesic segment from
$\overline W_a$, $\overline W_b$
to
$\overline W_v$
passes through
$\overline W_u$
and
$$d(\overline W_u,\overline W_v)<(c(k-1)+1/2)R.$$
We show that
$d(\overline W_a,\overline W_u)\geq R/2$.
Indeed, suppose that it is not true. Then there are
$x_1\in\overline W_a$, $x_2$, $z_2\in\overline W_u$,
$z_1\in\overline W_v$
such that
$|x_1x_2|< R/2$, $|z_1z_2|<(c(k-1)+1/2)R$.
Thus
$$|x_1z_1|\le|x_1x_2|+|x_2z_2|+|z_2z_1|<R/2+(c-1)R+(c(k-1)+1/2)R=ckR$$
in contradiction with
$\overline W_a\subset A_k$.

Similarly, we have
$d(\overline W_b,\overline W_u)\ge R/2$.
It follows from
$u\in ab\subset T$
that
$$d(\overline W_a,\overline W_b)
\ge d(\overline W_a,\overline W_u)+d(\overline W_u,\overline W_b)\ge R.$$

We restrict the constructed covering of
$\overline{W}$
to
$W\subset\overline{W}$,
i.e. for each member of the covering consider its intersection
with
$W$.
The new covering of
$W$
is again 2-colored,
$R$-disjoint
and
$c'R$-bounded
with
$c'=5c-2=88$.
We denote this covering by
$\Xx_W$, $\Xx_W=\Xx^1_W\cup \Xx^2_W$.

For each triple of numbers $(i,j,k)\in\{1,2\}^3$
consider the following 8-colored covering of
$X$,
$$\Xx:=\bigcup\{ \Xx^k_W \mid W\in\widetilde{\Ww}^{(i,j)},
\quad (i,j,k)\in \{1,2\}^3\}.$$

The covering
$\Xx$
is
$c'R$-bounded
because the families
$\Xx^k_W$, $k=1,2$,
are
$c'R$-bounded.
We prove that
$\Xx$
is
$R$-disjoint.
Indeed, let $U\in\Xx^k_{W_1},V\in \Xx^k_{W_2}$ be two
different sets of the same color for some triple of
numbers
$(i,j,k)\in \{1,2\}^3$,
where
$W_1,W_2\in \widetilde{\Ww}^{(i,j)}$.
If
$W_1\ne W_2$,
then
$d(U,V)\ge d(W_1,W_2)\ge R$,
since
$U\subset W_1$, $V\subset W_2$,
and the families
$\widetilde{\Ww}^{(i,j)}$
are
$R$-disjoint.
If
$W_1=W_2$,
then
$d(U,V)\ge R$
because
$\Xx^k_{W_1}$
is
$R$-disjoint.
This completes the proof of Theorem~\ref{t}.

\section{Applications}
Let us recall some definitions.
\begin{dfn}(Assouad \cite{As})
\textit{The  Assouad~-- Nagata dimension} of a metric
space $X$, $\andim X$, is the least integer number $n$
such that for each $R>0$ there exists an
$(n+1)$-colored, $R$-disjoint, $CR$-bounded covering of
the space $X$, where the number $C>0$ is independent
from $R$.
\end{dfn}

\begin{dfn}(Buyalo \cite{Buy2})
\textit{The linearly-controlled dimension} of a metric
space $X$, $\ldim X$, is the least integer number $n$,
such that for each $R\leq 1$ there exists an
$(n+1)$-colored, $R$-disjoint, $CR$-bounded covering of
the space $X$, where the number $C>0$ is independent
from $R$.
\end{dfn}

One easily sees that
$\andim X=\max\{\ldim X,\lasdim X\}$.

Recall that $\widetilde{M}$ is the universal cover of
the graph-manifold
$M$.
A remark to \cite[Proposition~2.7]{LS} implies that
$\ldim \widetilde{M}=\dim\widetilde{M}$,
so $3=\dim\widetilde{M}\leq \andim\widetilde{M}$.
Theorem~\ref{t} implies that
$\andim\widetilde{M}\leq 7$, which gives
Corollary~\ref{sle1}.
However, by Bell-Dranishnikov result (see \cite{BD})
$\asdim \widetilde{M}=3$,
which suggests following conjecture.

\begin{hyp}
Let $\widetilde{M}$ be the universal cover of the
graph-manifold. Then
$\andim\widetilde{M}=3$.
\end{hyp}

Corollary~\ref{sle:a} follows from
Corollary~\ref{sle1} and \cite[Theorem~1.3]{LS}.

\begin{dfn}
A metric space $X$ called \textit{Lipschitz
$n$-connected},
$n\ge 0$, if for each $m\in\{0,\ldots,n\}$
there exists
$\gamma>0$
such that every
$\lambda$-Lipschitz map
$f\colon S^m\to X$
permits a
$\gamma\lambda$-Lipschitz
extension
$\bar{f}\colon B^{m+1}\to X$,
where $S^m$
and
$B^{m+1}$
are the unit sphere and the unit disk in
$\R^{m+1}$
respectively.
\end{dfn}

It is proved in \cite{BH}, that Hadamard spaces are
Lipschitz $n$-connected for each $n\ge 0$. Thus the model space
$X$
is Lipschitz
$n$-connected
for each
$n\ge 0$. By Corollary \ref{sle1} the space
$X$
has finite Assouad-Nagata dimension. Corollary~\ref{sle2} now
follows from \cite[Theorem~1.5]{LS}.

\end{document}